\documentclass[10pt, article, titlepage,openany, rotate]{amsart}

\usepackage{amssymb, amsfonts,amscd,verbatim}
\usepackage{latexsym}
\usepackage{graphicx}
\usepackage{verbatim}   
\usepackage{color}      
\usepackage{subfigure}  

\usepackage{amsmath}    

\theoremstyle{plain}
\newtheorem{Prop}{Proposition}[section]
\newtheorem{Thm}[Prop]{Theorem}
\newtheorem{Cor}[Prop]{Corollary}

\newtheorem{Lem}[Prop]{Lemma}

\theoremstyle{definition}
\newtheorem{Def}[Prop]{Definition}

\theoremstyle{remark}

\newtheorem{Problem}[Prop]{\bf Problem}

\newcommand{\map}{\operatorname{Map}\nolimits}

\def\dim{\mathop{\roman{dim}}}
\def\int{\mathop{\roman{int}}}

\def\1{^{-1}}

\def\dim{\text{dim}}

\def\Ext{\text{Ext}}
\def\Tor{\text{Tor}}
\def\Hom{\text{Hom}}

\def\Ab{\text{Ab}}

\def\PPP{{\mathcal P}}

\newcommand{\ZZ}{\mathbb Z}
\newcommand{\QQ}{\mathbb Q}

\def\dokaz{{\bf Proof. }}
\numberwithin{equation}{section}
\begin{document}
\title[
Hurewicz-Serre Theorem in extension theory
]%
   {Hurewicz-Serre Theorem in extension theory}

\author{M.~Cencelj}
\address{Fakulteta za Matematiko in Fiziko,
Univerza v Ljubljani,
Jadranska ulica 19,
SI-1111 Ljubljana,
Slovenija }
\email{matija.cencelj@guest.arnes.si}

\author{J.~Dydak}
\address{University of Tennessee, Knoxville, TN 37996, USA}
\email{dydak@math.utk.edu}

\author{A.~Mitra}
\address{University of Tennessee, Knoxville, TN 37996, USA}
\email{ajmitra@math.utk.edu}

\author{A.~Vavpeti\v c}
\address{Fakulteta za Matematiko in Fiziko,
Univerza v Ljubljani,
Jadranska ulica 19,
SI-1111 Ljubljana,
Slovenija }
\email{ales.vavpetic@fmf.uni-lj.si}

\date{ January 21, 2006}
\keywords{Extension dimension, cohomological dimension, absolute extensor,
nilpotent groups}

\subjclass{ Primary: 54F45; Secondary: 55M10, 54C65
}

\thanks{Supported in part by the Slovenian-USA research grant BI--US/05-06/002 and the ARRS
research project No. J1--6128--0101--04. The second-named author was partially supported
by Grant No.2004047  from the United States-Israel Binational Science
Foundation (BSF),  Jerusalem, Israel}

\begin{abstract}
The paper is devoted to generalizations of Cencelj-Dranishnikov theorems
relating extension properties of nilpotent CW complexes to its homology
groups.
Here are the main results of the paper:
\begin{Thm}
Suppose $L$ is a nilpotent CW complex and $F$
is the homotopy fiber of the inclusion $i$ of $L$ into
its infinite symmetric product $SP(L)$.
If $X$ is a metrizable space such that $X\tau K(H_k(L),k)$ for all $k\ge 1$,
then $X\tau K(\pi_k(F),k)$ and $X\tau K(\pi_k(L),k)$ for all $k\ge 2$.
\end{Thm}

\begin{Thm}
Let $X$ be a metrizable space such that $\dim(X) < \infty$ or $X\in ANR$.
Suppose $L$ is a nilpotent CW complex and $SP(L)$ is its infinite symmetric product.
If $X\tau SP(L)$,
then $X\tau L$ in the following cases:
\begin{itemize}
\item[a.] $H_1(L)$ is  finitely generated.
\item[b.] $H_1(L)$ is a torsion group.
\end{itemize}
\end{Thm}

\end{abstract}

\maketitle

\medskip
\medskip
\tableofcontents

\section{Introduction}\label{SectionIntro}

Recall that $X\tau L$ and $L\in AE(X)$ are shortcuts to the statement: $L$ is an absolute extensor
of $X$. The geometric result being repeatedly used in the paper is the consequence of Theorem 1.9 in \cite{DrDy1}: if $F$ is the homotopy fiber of $L\to K$ and $X\tau F$,
then $X\tau L$ is equivalent to $X\tau K$.

Given a metrizable space $X$ and a connected CW complex $L$
consider the following conditions:
\begin{enumerate}
\item $X\tau L$.
\item $X\tau K(\pi_n(L),n)$ for all $n\ge 1$.
\item $X\tau K(H_n(L),n)$ for all $n\ge 1$.
\end{enumerate}

It is well-known that (1) implies (3) as proved by Dranishnikov \cite{DR0} for
$X$ compact and for arbitrary $X$ in [\cite{Dy2}, Theorem 3.4].
The difficulty in generalizing results in cohomological dimension theory
from compact spaces to arbitrary metrizable spaces
usually lies in the fact that the First Bockstein Theorem does not hold for
metric spaces.

By a Hurewicz-Serre Theorem in Extension Theory we mean any
result showing (3) implies (2). However, in practice we are really
interested in arriving at (1).

Here is the main problem we are interested in:
\begin{Problem}\label{MainProblemForNilpotentCWComplexes} 
Suppose $X$ is a metrizable space such that $X\tau K(H_n(L),n)$ for all $n\ge 1$.
If $L$ is nilpotent, does  $X\tau K(\pi_n(L),n)$ hold for all $n\ge 1$?
\end{Problem}

Even a specialized version of \ref{MainProblemForNilpotentCWComplexes} 
is open:

\begin{Problem}\label{IsPiOneAbsExtensorProblem}
Suppose $L$ is a nilpotent CW complex.
If $X$ is a metrizable space such that $X\tau K(H_k(L),k)$ for all $k\ge 1$,
does $X\tau K(\pi_1(L),1)$ hold?
\end{Problem}

Notice that it is not sufficient to assume $X\tau K(H_1(L),1)$
in \ref{IsPiOneAbsExtensorProblem}. Namely, take the group $G$ from \cite{War}
whose abelianization is $\QQ\oplus \QQ$ and whose commutator group
is $\ZZ/p^\infty$. Pick a compactum $X$ so that $\dim_\QQ(X)=1$
and $\dim_{\ZZ/p^\infty}(X)=2$.
The complex $L=K(G,1)$ is nilpotent and $X\tau K(H_1(L),1)$
but $X\tau L$ does not hold. Indeed, as $\pi_1(L)$ is not $\bar p$-local
and $H_1(L;\ZZ/p^\infty)=0$, \ref{PiOneIsPLocal} says $H_2(L;\ZZ/p^\infty)\ne 0$
which means $H_2(L)/\Tor(H_2(L))$ is not $p$-divisible.
If $X\tau L$, then $X\tau K(H_2(L),2)$ and $\dim_{\ZZ_{(p)}}(X)\leq 2$
and that, in combination with $\dim_\QQ(X)=1$, implies $\dim_{\ZZ/p^\infty}(X)\leq 1$,
a contradiction. See more in \cite{CDMV} about Bockstein First Theorem for nilpotent groups.

However, if $H_1(L)$ is a torsion group, then 
the answer to \ref{IsPiOneAbsExtensorProblem}
is positive. 

\begin{Lem}\label{ExtensionWRTFundGroupOfComplex}
Suppose $N$ is a nilpotent group.
If $X\tau K(\Ab(N),1)$ for some metrizable space $X$ and
$\Ab(N)$ is a torsion group, then $X\tau K(N,1)$.
\end{Lem}
\dokaz 
We will prove \ref{ExtensionWRTFundGroupOfComplex} by induction on the nilpotency class $n$ of $N$.
Let $\Gamma^nN=\Gamma^n$. Notice $N/\Gamma^n$ is a nilpotent group of class $(n-1)$
whose abelianization is an image of $\Ab(N)$. Thus $X\tau K(N/\Gamma^n,1)$. The epimorphism
$$
\otimes^n {\mathrm{Ab}}N\longrightarrow \Gamma^nN=\Gamma^n
$$
implies $X\tau K(\Gamma^n,1)$, so the fact that $N$ is a central extension
$$
1\rightarrow \Gamma^n\rightarrow N\rightarrow N/\Gamma^n\rightarrow 1
$$
concludes the proof.
\hfill $\blacksquare$

For the sake of completeness let us show (2) is always stronger than (3).

\begin{Prop}\label{HomotopyImpliesHomology} 
Suppose $X$ is a metrizable space 
and $L$ is a connected CW complex. If  $X\tau K(\pi_n(L),n)$ hold for all $n\ge 1$,
then $X\tau K(H_n(L),n)$ for all $n\ge 1$.
\end{Prop}
\dokaz Let $L_n$ be the CW complex obtained from $L$ by killing all
homotopy groups higher than $n$. Since $L_n$ is obtained from $L$ by attaching
$k$-cells for $k> n+1$, $H_n(L_n)=H_n(L)$.
Also, one has $X\tau L_n$ as $X\tau K(\pi_i(L_n),i)$ holds for all $i$
and only finitely many homotopy groups of $L_n$ are non-trivial
(see Theorem G of \cite{Dy1}). Therefore $X\tau K(H_n(L_n),n)$.
\hfill $\blacksquare$

Also, (1) is always stronger than (2) provided the issue of the fundamental group is avoided.

\begin{Prop}\label{XtauLImpliesHomotopy} 
Suppose $X$ is a metrizable space 
and $L$ is a connected CW complex. If  $X\tau K(\pi_1(L),1)$ and $X\tau L$,
then $X\tau K(\pi_n(L),n)$ hold for all $n\ge 1$.
\end{Prop}
\dokaz Notice that the homotopy fiber of the covering
projection $\tilde L\to L$ is $K(\pi_1(L),1)$. Therefore $X\tau \tilde L$
and $X\tau K(H_n(\tilde L),n)$ for all $n$.
By Theorem F of \cite{Dy1} (see also \ref{ExtensionWRTHomotopyOfFiber} of this paper) one has
$X\tau K(\pi_n(\tilde L),n)=K(\pi_n(L),n)$ for all $n\ge 2$.
\hfill $\blacksquare$

\begin{Def}\label{KnoxvilleSpace}
$X$ is called a {\it Knoxville space} if it is metrizable and for any connected CW complex $L$
the conditions $X\tau K(\pi_n(L),n)$ for all $n\ge 1$ imply $X\tau L$.
\end{Def}

\begin{Problem}\label{CharKnoxvilleSpacesProblem} 
Characterize Knoxville spaces.
\end{Problem}

It follows from Theorem G of \cite{Dy1} that any finitely dimensional $X$
or any $X\in ANR$ is a Knoxville space. Also, it is easy to see that any countable
union of closed Knoxville subspaces is a Knoxville subspace.

\section{Properties of the homotopy fiber of $L\to SP(L)$}

Notice that condition (3) of Section~\ref{SectionIntro}
is equivalent to $X\tau SP(L)$ as $SP(L)$ is the weak product
of $K(H_n(L),n)$ for all $n\ge 1$ according to the famous theorem of Dold and Thom \cite{D-T}.
Since we are interested in deriving $X\tau L$ it makes sense to ponder
the stronger condition $X\tau F$, where $F$ is the homotopy fiber
of the inclusion $L\to SP(L)$. That is the main idea of the whole paper
and in this section we concentrate on basic properties of $F$
and its homotopy groups.

\begin{Prop}\label{FiberIsNilpotent}
Suppose $L$ is a CW complex and $F$
is the homotopy fiber of the inclusion $i$ of $L$ into
its infinite symmetric product $SP(L)$.
If $L$ is nilpotent, then $F$ is nilpotent.
\end{Prop}
\dokaz
The homotopy sequence 
$$
\cdots \to \pi_{n}(F)\stackrel{j_*}{\to} \pi_{n}(L)\stackrel{i_*}{\to}\pi_{n}(SP(L))\stackrel\partial{\to} \pi_{n-1}(F)\to\cdots
$$
of the fibration $F\stackrel{j}{\to} L\stackrel{i}{\to} SP(L)$ is sequence
of $\pi_1(L)$-modules \cite[Proposition 8$^{bis}$.2]{Mac}.
For the action of $\pi_1(L)$ on $\pi_n(F)$ which is described
in the proof of \cite[Proposition 8$^{bis}$.2]{Mac} holds
$g\cdot \alpha =j_*(g)\cdot\alpha$ for $\alpha\in\pi_n(F)$ and $g\in\pi_1(F)$.

Let $I_F$ and $I_L$ be the augmentation ideals of group rings
$\ZZ[\pi_1(F)]$ and $\ZZ[\pi_1(L)]$, respectively.
Because $L$ is a nilpotent, there is an integer $c$, such that
$(I_L)^c\pi_n(L)=0$.
Let $\eta\in (I_F)^c$ and $\alpha\in\pi_n(F)$.
Then $j_*(\eta\alpha)=j_*(\eta) j_*(\alpha)=0$, because $j_*(\eta)\in (I_L)^c$.
Thus there exists $\beta\in\pi_{n+1}(SP(L))$, such that 
$\partial\beta=\eta\alpha$.
Let $g\in\pi_1(F)$. Then $(j_*(g)-1)\in I_L$ and
$$
\partial((j_*(g)-1)\beta)=(j_*(g)-1)\partial\beta=(j_*(g)-1)\eta\alpha=
(g-1)\eta\alpha.
$$
The action of $\pi_1(L)$ on $\pi_n(SP(L))$ is defined as
$l\gamma=i_*(l)\gamma$ for $l\in \pi_1(L)$ and $\gamma\in \pi_n(SP(L))$.
Hence 
$$
(j_*(g)-1)\beta=(i_*j_*(g)-1)\beta=(1-1)\beta=0,
$$
therefore $(g-1)\eta\alpha=0$. This shows that
$(I_F)^{c+1} \pi_n(F)=0$, so the space $F$ is nilpotent.
\hfill $\blacksquare$

\begin{Prop}\label{SerreThmForFiber}
Suppose $L$ is a nilpotent CW complex and $F$
is the homotopy fiber of the inclusion $i$ of $L$ into
its infinite symmetric product $SP(L)$.
If $\PPP$ is a set of primes such that $H_k(L)$ is a $\PPP$-torsion group
for all $k\leq n$, where $n\ge 1$ is given,
then $\pi_k(F)$ is a $\PPP$-torsion group for all $1\leq k\leq n+1$.
\end{Prop}
\dokaz Let $\PPP'$ be the complement of $\PPP$ in the set of all primes.
Consider the localization $L_{(\PPP')}$ of $L$ at $\PPP'$.
Notice that $L_{(\PPP')}$ is $n$-connected, so the Hurewicz homomorphism
$\phi_k:\pi_k(L_{(\PPP')})\to H_k(L_{(\PPP')})$ is an isomorphism for $k=n+1$
and an epimorphism for $k=n+2$.
Let us split the exact sequence $\ldots\to \pi_k(F)\to \pi_k(L)\to H_k(L)\to\ldots$
into $\ldots\to \pi_2(F)\to \pi_2(L)\to H_2(L)\to A\to 0$ and $1\to A\to \pi_1(F)\to B\to 1$,
where $B$ is the commutator subgroup of $\pi_1(L)$.
Localizing the first sequence at $\PPP'$ yields $A$ being a $\PPP$-torsion group
and $\pi_k(F)$ being $\PPP$-torsion for $2\leq k\leq n+1$. Since $B$ is $\PPP$-torsion,
\ref{SerreThmForFiber} follows.
\hfill $\blacksquare$

\begin{Cor}\label{HurSerreThmInTorsionCase}
Suppose $L$ is a nilpotent CW complex and $F$
is the homotopy fiber of the inclusion $i$ of $L$ into
its infinite symmetric product $SP(L)$.
If $n > 1$ is a number such that $H_k(L)$ is a torsion group for all $k < n$,
then for any metrizable space $X$ the conditions $X\tau K(H_k(L),k)$
for all $k\leq n$ 
imply $X\tau K(\pi_k(F),k)$ for all $1\leq k\leq n$.
\end{Cor}
\dokaz The case $k=1$ is taken care of by \ref{ExtensionWRTFundGroupOfComplex}.
If $p$-torsion of $\pi_k(F)$ is not trivial,
then \ref{SerreThmForFiber} implies that $p$-torsion of $H_m(L)$
is not trivial for some $m < k$. Therefore
$X\tau K(\ZZ/p^\infty,m)$ and $X\tau K(\ZZ/p,m+1)$.
This implies $X\tau K(G,k)$ for all $G$ in the Bockstein basis of $\pi_k(F)$
resulting in $X\tau K(\pi_k(F),k)$.
\hfill $\blacksquare$

\section{Homotopy groups with coefficients}
Given a countable Abelian group $G$ consider a pointed compactum $P_n(G)$
such that its integral cohomology is concentrated in dimension $n$
and equals $G$. The $n$-th homotopy group $\pi_n(L;G)$
of a pointed CW complex $L$ is defined in \cite{Nei} to be the set
$[P_n(G),L]$ of pointed homotopy classes from $P_n(G)$ to $L$.
If $P_{n-1}(G)$ exists (that is always true if $n > 2$ or $G$ is torsion free and $n\ge 2$), then $P_n(G)$ could be taken as the suspension $\Sigma P_{n-1}(G)$
of $P_{n-1}(G)$
with the resulting group structure on $\pi_n(L,G)$.

If one puts $D=P_2(G)$ (or $D=P_1(G)$ if $G$ is torsion-free), then one can analyze homotopy groups of $L^D=\map(D,L)$ and realize that $\pi_n(L^D)=\pi_{n+2}(L;G)$ (respectively, $\pi_n(L^D)=\pi_{n+1}(L;G)$). 
Therefore, given a Hurewicz fibration $F\to E\to B$, one concludes there is a long exact sequence
$\ldots\to \pi_n(F;G)\to \pi_n(E;G)\to\pi_n(G;G)\to \pi_{n-1}(F;G)\to\ldots$ (see \cite{Nei} for the special case of $G=\ZZ/m$) because $F^D\to E^D\to B^D$ is a Serre fibration.

\par In the special case of $G=\ZZ/m$ one can pick the Moore space $D=M(\ZZ/m,1)$ for $P_2(G)$.
In that case one has a Serre fibration (that follows from the Homotopy Extension Theorem) $\map(S^2,L)\to \map(D,L)\to \map(S^1,L)$ where $S^1$ is the 1-skeleton of $D$ and $S^2=D/S^1$. The map $\map(D,L)\to \map(S^1,L)$ is simply restriction induced. Since the boundary homomorphism
$\pi_{n+1}(B)\to \pi_n(F)$ in that case amounts to multiplication by $m$ from
$\pi_{n+1}(\map(S^1,L))=\pi_{n+2}(L)$ to $\pi_{n}(\map(S^2,L))=\pi_{n+2}(L)$, one concludes the following
(see \cite{Nei} for another way of deriving an equivalent result):

\begin{Lem} \label{UnivCoeffForHomotopyGroupsModp}
Let $D=M(\ZZ/m,1)$ for some $m\ge 2$.
For each pointed CW complex $L$ and each $n\ge 0$ one has a natural exact sequence
$$0\to \pi_{n+2}(L)\otimes \ZZ/m\to \pi_n(L^D)\to  \pi_{n+1}(L)\ast \ZZ/m\to 0,$$
where $ \pi_{1}(L)\ast \ZZ/m$ is $\{x\in  \pi_{1}(L)| x^m=1\}$.
\end{Lem}

We are interested in homotopy groups with coefficients in $\ZZ/p^\infty$,
the direct limit of $\ZZ/p\to\ZZ/p^2\to\ldots$. Notice that one can construct
$P_2(\ZZ/p^\infty)$ as the inverse limit of $\ldots\to M(\ZZ/p^{n+1},1)\to M(\ZZ/p^n,1)\to\ldots\to M(\ZZ/p,1)$
which can be viewed as $M(\hat Z_p,1)$, the Moore space for the $p$-adic integers $\hat Z_p$
in terms of Steenrod homology.
In that case \ref{UnivCoeffForHomotopyGroupsModp}
becomes
\begin{Lem} \label{UnivCoeffForHomotopyGroupsModpInfty}
Let $p$ be prime.
For each pointed CW complex $L$ and each $n\ge 0$ one has a natural exact sequence
$$0\to \pi_{n+2}(L)\otimes \ZZ/p^\infty\to \pi_{n+2}(L;\ZZ/p^\infty)\to  \pi_{n+1}(L)\ast \ZZ/p^\infty\to 0,$$
where $ \pi_{1}(L)\ast \ZZ/p^\infty$ is $\{x\in  \pi_{1}(L)| x^{p^k}=1\text{ for some }k\ge 1\}$.
\end{Lem}

As a consequence of \ref{UnivCoeffForHomotopyGroupsModp},
\ref{UnivCoeffForHomotopyGroupsModpInfty}, and Dold-Thom Theorem \cite{D-T}
($\pi_n(SP(L)=H_n(L)$)
one can get that $\pi_n(SP(L);G)=H_n(L;G)$
for all $n$ and $G=\ZZ/p$ or $G=\ZZ/p^\infty$.

\begin{Prop} \label{HomotopyOfTheFiberModp}
Suppose $L$ is a nilpotent CW complex whose fundamental group is $\bar p$-local
for some prime $p$. Let $F$ be the homotopy fiber
of the inclusion $i:L\to SP(L)$ of $L$ into its infinite symmetric product.
If $H_k(L,\ZZ/p)=0$ for $k \leq n$, where $n\ge 1$, then
$\pi_{k}(F,\ZZ/p)=0$ for $2\leq k\leq n+1$ and $\pi_{k}(L,\ZZ/p)=0$ for $2\leq k\leq n$.
\end{Prop}
\dokaz  Let $\tilde L$ be the universal cover of $L$
and let $\pi:\tilde L\to L$ be the covering projection. 
Recall that every nilpotent CW complex $L$ has the $p$-completion
$L_p$ with the map $L\to L_p$ inducing isomorphisms of all
homology groups with coefficients in $\ZZ/p$
such that one has a natural exact sequence
$$0\to \Ext(\ZZ/p^\infty,\pi_n(L))\to \pi_n(L_p)\to \Hom(\ZZ/p^\infty,\pi_{n-1}(L))\to 0$$
for all $n\ge 1$ (see [\cite{GoerssJardine}, Theorem 3.7 on p.416]).
By \ref{ExtIsZero} and \ref{HomIsZero} one has
$\Ext(\ZZ/p^\infty,\pi_1(L))= \Hom(\ZZ/p^\infty,\pi_{1}(L))=0$, so
 the induced map
$\hat{\tilde L}_p\to \hat L_p$ is a homotopy equivalence. 
Therefore it induces isomorphism of homology mod $p$ and $\pi$ induces
isomorphism of homology mod $p$. However, $H_2(\tilde L;\ZZ/p)=\pi_2(\tilde L;\ZZ/p)=\pi_2(L;\ZZ/p)$
and $\pi_3(L)\to H_3(\tilde L)$ is an epimorphism resulting
in $\pi_3(L;\ZZ/p)\to H_3(L;\ZZ/p)$ being an epimorphism.
By exactness of mod $p$ groups of a fibration one gets $\pi_{2}(F,\ZZ/p)=0$.
That proves \ref{HomotopyOfTheFiberModp} for $n=1$.
\par If $n > 1$ we apply mod $p$ Hurewicz Theorem of \cite{Nei} to $\tilde L$
to conclude $\pi_{n+1}(\tilde L;\ZZ/p)\to H_{n+1}(\tilde L;\ZZ/p)$ is an isomorphism
and $\pi_{n+2}(\tilde L;\ZZ/p)\to H_{n+2}(\tilde L;\ZZ/p)$ is an epimorphism.
Consequently, $\pi_{n+1}(L;\ZZ/p)\to H_{n+1}(L;\ZZ/p)$ is an isomorphism
and $\pi_{n+2}(L;\ZZ/p)\to H_{n+2}(L;\ZZ/p)$ is an epimorphism. Thus $\pi_{n+1}(F;\ZZ/p)=0$.
\hfill $\blacksquare$

\begin{Cor} \label{HomotopyOfTheFiberModpInfty}
Suppose $L$ is a nilpotent CW complex whose fundamental group is $\bar p$-local
for some prime $p$. Let $F$ be the homotopy fiber
of the inclusion $i:L\to SP(L)$ of $L$ into its infinite symmetric product.
If $H_k(L,\ZZ/p^\infty)=0$ for $k \leq n$, where $n\ge 2$, then
$\pi_{k}(L,\ZZ/p^\infty)=\pi_{k}(F,\ZZ/p^\infty)=0$ for $2\leq k\leq n$.
\end{Cor}
\dokaz Case 1: $n > 2$. Notice $H_k(L,\ZZ/p)=0$ for $k \leq n-1$ resulting in
$\pi_{k}(F,\ZZ/p)=0$ for $2\leq k\leq n$. Hence
$\pi_{k}(F,\ZZ/p^\infty)=0$ for $2\leq k\leq n$
and from an exact sequence we get $\pi_{k}(L,\ZZ/p^\infty)=0$ for $2\leq k\leq n$.
\par Case 2: $n=2$.
Let $\tilde L$ be the universal cover of $L$
and let $\pi:\tilde L\to L$ be the covering projection. Notice that the induced map
$\hat{\tilde L}_p\to \hat L_p$ is a homotopy equivalence. 
Therefore it induces isomorphism of homology mod $p$ and $\pi$ induces
isomorphism of homology mod $p$. However, $H_2(\tilde L;\ZZ/p^\infty)=\pi_2(L;\ZZ/p^\infty)$
and $\pi_3(L)\to H_3(\tilde L)$ is an epimorphism resulting
in $\pi_3(L;\ZZ/p^\infty)\to H_3(L;\ZZ/p^\infty)$ being an epimorphism.
By exactness of mod $p$ groups of a fibration one gets $\pi_{2}(F,\ZZ/p^\infty)=\pi_{2}(L,\ZZ/p^\infty)=0$.
\hfill $\blacksquare$

\begin{Prop} \label{HomotopyOfTheFiberModZLocalized}
Suppose $L$ is a nilpotent CW complex, $F$ is the homotopy fiber
of the inclusion $i:L\to SP(L)$ of $L$ into its infinite symmetric product,
and $n\ge 2$.
If $H_k(L,\ZZ_{(\PPP)})=0$ for $k \leq n$, then
$\pi_{k}(F,\ZZ_{(\PPP)})=0$ for $2\leq k\leq n+1$
and $\pi_{k}(L,\ZZ_{(\PPP)})=0$ for $2\leq k\leq n$.
\end{Prop}
\dokaz Let $\PPP'$ be the complement of $\PPP$ in the set of all primes.
If $H_k(L,\ZZ_{(\PPP)})=0$ for $k \leq n$, then $H_1(L)$ is a $\PPP'$-torsion group
resulting in $\pi_1(L)$ being a $\PPP'$-torsion group.
Let $L_{(\PPP)}$ be the $\PPP$-localization of $L$. It is $n$-connected,
so by the classical Hurewicz Theorem $\pi_k(L_{(\PPP)})\to H_k(L_{(\PPP)})$
is an isomorphism for $k\leq n+1$ and an epimorphism for $k=n+2$.
That corresponds to $\pi_k(L;\ZZ_{(\PPP)})\to H_k(L;\ZZ_{(\PPP)})$
being an isomorphism for $k\leq n+1$ and an epimorphism for $k=n+2$.
In view of exactness of $\ldots \to \pi_k(F;\ZZ_{(\PPP)})\to \pi_k(L;\ZZ_{(\PPP)})\to H_k(L;\ZZ_{(\PPP)})\to\ldots$,  \ref{HomotopyOfTheFiberModZLocalized} follows.
\hfill $\blacksquare$

\section{Main results}

\begin{Lem} \label{CasesWhenCohDimBig}
Suppose $X$ is a metrizable space, $G$ is an Abelian group, and $n\ge 1$.
If $\dim_G(X) > n$, then one of the following conditions holds:
\begin{itemize}
\item[a.] $\dim_\QQ(X)\ge n+1$ and $G$ is not a torsion group.
\item[b.] There is a prime $p$ such that $G\otimes \ZZ/p^\infty\ne 0$
and $\dim_\QQ(X)\leq n$, $\dim_{\ZZ/p^\infty}(X)\ge n$.
\item[c.] There is a prime $p$ such that $G$ is $p$-divisible, $\Tor_p(G)\ne 0$, and $\dim_{\ZZ/p^\infty}(X)\ge n+1$.
\item[d.] There is a prime $p$ such that $\Tor_p(G)$ is not $p$-divisible and $\dim_{\ZZ/p}(X)\ge n+1$.
\end{itemize}
\end{Lem}
\dokaz Let $\dim_G(X)=m$. Suppose none of (a)-(d) holds. 
According to Part (b) of Theorem B of \cite{Dy1} one has
$m=\dim_{G/\Tor(G)}(X)$ or $m=\dim_{\Tor(G)}(X)$.
If $\dim_{\Tor(G)}(X)\ge n+1$, then, according to Part (a) of Theorem B of \cite{Dy1},
there is a prime $p$ such that either $\Tor(G)$ is $p$-divisible,
$\Tor_p(G)\ne 0$, and $\dim_{\Tor(G)}(X)=\dim_{\ZZ/p^\infty}(X)$
(in which case (c) holds) or $\Tor(G)$ is not $p$-divisible,
$\Tor_p(G)\ne 0$, and $\dim_{\Tor(G)}(X)=\dim_{\ZZ/p}(X)$
(in which case (d) holds).
Therefore $m=\dim_{G/\Tor(G)}(X)$  and $\dim_{\Tor(G)}(X)\leq n$.
In particular $G$ is not a torsion group, so $\dim_\QQ(X)\leq n$ as (a)
fails to hold.

\par Consider $\PPP=\{p | G\otimes \ZZ/p^\infty\ne 0\}$,
the set of primes $p$ such that $G/\Tor(G)$ is not $p$-divisible.
It is shown in \cite{Dy1} (Part (f) of Theorem B) that $\dim_{\ZZ_{(\PPP)}}(X)\ge \dim_{G/\Tor(G)}(X)$,
so $\dim_{\ZZ_{(\PPP)}}(X)\ge m$.
As (b) does not hold, one has $\dim_{\ZZ/p^\infty}(X)\leq n-1$ for all $p\in \PPP$.
From the exact sequence $$0\to \ZZ_{(\PPP)}\to \QQ\to \bigoplus\limits_{p\in\PPP}\ZZ/p^\infty\to 0$$
one concludes that the homotopy fiber of $K(\ZZ_{(\PPP)},m-1)\to K( \QQ, m-1)$
is $K( \bigoplus\limits_{p\in\PPP}\ZZ/p^\infty,m-2)$.
Since $m-2\ge n-1$, $X\tau K( \bigoplus\limits_{p\in\PPP}\ZZ/p^\infty,m-2)$
which implies $X\tau K(\ZZ_{(\PPP)},m-1)$
as $X\tau K(\QQ,m-1)$ is true. Thus
$\dim_{\ZZ_{(\PPP)}}(X) \leq m-1$, a contradiction.
\hfill $\blacksquare$

\begin{Thm}\label{ExtensionWRTHomotopyOfFiber}
Suppose $L$ is a nilpotent CW complex and $F$
is the homotopy fiber of the inclusion $i$ of $L$ into
its infinite symmetric product $SP(L)$.
If $X$ is a metrizable space such that $X\tau K(H_k(L),k)$ for all $k\ge 1$,
then $X\tau K(\pi_k(F),k)$ and $X\tau K(\pi_k(L),k)$ for all $k\ge 2$.
\end{Thm}
\dokaz Suppose $n\ge 2$ is the smallest natural number such that 
$X\tau K(\pi_k(F),k)$ fails (similar argument in case $X\tau K(\pi_k(L),k)$ fails). 
By \ref{CasesWhenCohDimBig} one of the following cases holds 
for $G=\pi_n(F)$:
\begin{itemize}
\item[a.] $\dim_\QQ(X)\ge n+1$ and $G$ is not a torsion group.
\item[b.] There is a prime $p$ such that $G\otimes \ZZ/p^\infty\ne 0$
and $\dim_\QQ(X)\leq n$, $\dim_{\ZZ/p^\infty}(X)\ge n$.
\item[c.] There is a prime $p$ such that $G$ is $p$-divisible, $\Tor_p(G)\ne 0$, and $\dim_{\ZZ/p^\infty}(X)\ge n+1$.
\item[d.] There is a prime $p$ such that $\Tor_p(G)$ is not $p$-divisible and $\dim_{\ZZ/p}(X)\ge n+1$.
\end{itemize}

\par {\bf Case 1}: $\dim_{\QQ}(X)\leq n-1$. Now only (b)-(d) are possible. Let $p$ be the prime from one
of those cases. Notice $H_k(L)$ is $\bar p$-local for $k\leq n-1$
as otherwise $\dim_{\ZZ/p^\infty}(X)\leq n-1$ and $\dim_{\ZZ/p}(X)\leq n$
so none of (b)-(d) would be valid.
Another observation is $H_n(L)\otimes \ZZ/p^\infty=0$. Indeed, 
$H_n(L)\otimes \ZZ/p^\infty\ne 0$ leads to $H_n(L)/\Tor(H_n(L))$ not being $p$-divisible
in which case [\cite{Dy1}, Part (d) of Theorem B] implies $\dim_{\hat \ZZ_p}(X)\leq n$
as $\dim_{H_n(L)}(X)\leq n$. Therefore $\dim_{\ZZ_{(p)}}(X)\leq n$
(see Part (e) of Theorem B in \cite{Dy1}) and
$\dim_{\ZZ/p^\infty}(X)\leq  \max(\dim_\QQ(X),\dim_{\ZZ_{(p)}}(X)-1)\leq n-1$,
a contradiction.
\par $\pi_1(L)$ is $\bar p$-local by \ref{PiOneIsPLocal} and $G\otimes \ZZ/p^\infty=0$ by  \ref{HomotopyOfTheFiberModpInfty}.
That means (b) is not possible. If $H_n(L)$ is not $p$-divisible, then $\dim_{\ZZ/p}(X)\leq n$
and neither (c) nor (d) would be possible. Thus $H_k(L;\ZZ/p)=0$ for $k\leq n$
resulting in $G$ being $p$-divisible by \ref{HomotopyOfTheFiberModp}.
That means only (c) is possible. In addition, $\Tor_p(H_n(L))=0$.
Also $H_{n+1}(L)\otimes \ZZ/p^\infty=0$ (otherwise 
$\dim_{\ZZ_{(p)}}(X)\leq n+1$ and
$\dim_{\ZZ/p^\infty}(X)\leq  \max(\dim_\QQ(X),\dim_{\ZZ_{(p)}}(X)-1)\leq n$). Thus $H_k(L;\ZZ/p^\infty)=0$ for $k\leq n+1$.
By \ref{HomotopyOfTheFiberModpInfty} $\Tor_p(G)=0$, a contradiction.
\par {\bf Case 2}: $\dim_{\QQ}(X) > n-1$. 
By \ref{SerreThmForFiber} the group $G$ is $\PPP$-torsion such that
$\dim_{\ZZ/p^\infty}(X)\leq n-1$ for all $p\in\PPP$ which implies $\dim_G(X)\leq n$, a contradiction.
\hfill $\blacksquare$

\begin{Cor}
Suppose $L$ is a nilpotent CW complex such that $\pi_n(L)=\pi_{n+1}(L)=0$
for some $n\ge 1$.
If $X\tau SP(L)$ for some metrizable space $X$, then $X\tau K(H_{n+1}(L),n)$.
\end{Cor}
\dokaz If $n=1$, then $H_{n+1}(L)=0$, so assume $n\ge 2$. Notice that
$\pi_n(F)=H_{n+1}(L)$, where $F$ is the homotopy fiber of $i:L\to SP(L)$.
\hfill $\blacksquare$

\begin{Lem}\label{ExtensionWRTFundGroupOfFiber}
Suppose $L$ is a nilpotent CW complex and $F$
is the homotopy fiber of the inclusion $i$ of $L$ into
its infinite symmetric product $SP(L)$.
If $X\tau K(H_1(L),1)$ for some metrizable space $X$ and $H_1(L)$ is  finitely generated,
then $X\tau K(\pi_1(F),1)$
\end{Lem}
\dokaz If $H_1(L)$ is finitely generated and non-torsion,
then $X$ is at most $1$-dimensional in which case $X\tau L$ for all connected CW complexes. Therefore assume $H_1(L)$ is a torsion group and (see \ref{SerreThmForFiber}) there is an exact sequence $1\to A\to \pi_1(F)\to B\to 1$
such that $A$ and $B$ are $\PPP$-torsion groups, where $\PPP=\{p | \Tor_p(H_1(L))\ne 0\}$.
Notice $H_1(L)$ does not contain $\ZZ/p^\infty$ for any prime $p$, so $X\tau K(A,1)$ and $X\tau K(B,1)$ which implies $X\tau K(\pi_1(F),1)$.
\hfill $\blacksquare$

\begin{Thm}
Let $X$ be a metrizable space such that $\dim(X) < \infty$ or $X\in ANR$.
Suppose $L$ is a nilpotent CW complex and $SP(L)$ is its infinite symmetric product.
If $X\tau SP(L)$,
then $X\tau L$ in the following cases:
\begin{itemize}
\item[a.] $H_1(L)$ is  finitely generated.
\item[b.] $H_1(L)$ is a torsion group.
\end{itemize}
\end{Thm} 
\dokaz a. By \ref{ExtensionWRTHomotopyOfFiber} and \ref{ExtensionWRTFundGroupOfFiber}
one concludes $X\tau K(\pi_n(F),n)$ for all $n\ge 1$.
Theorem G of \cite{Dy1} gives $X\tau F$ which implies $X\tau L$.
\par b.  By \ref{ExtensionWRTHomotopyOfFiber} and \ref{ExtensionWRTFundGroupOfComplex}
one concludes $X\tau K(\pi_n(L),n)$ for all $n\ge 1$.
Theorem G of \cite{Dy1} yields $X\tau  L$.
\hfill $\blacksquare$

\section{Appendix}
\begin{Lem} \label{AbDivisibleImpliesGIsDivisible}
Suppose $N$ is a nilpotent group and $p$ is a prime.
$\Ab(N)$ is $p$-divisible if and only if $N$ is $p$-divisible.
\end{Lem}
\dokaz
If $N$ is $p$-divisible clearly so is its abelianization. 
We will prove the converse by induction on the nilpotency class $n$ of $N$.
Let $\Gamma^nN=\Gamma^n$. Notice $N/\Gamma^n$ is a nilpotent group of class $(n-1)$
whose abelianization is $p$-divisible. Thus it is $p$-divisible. The epimorphism
$$
\otimes^n {\mathrm{Ab}}N\longrightarrow \Gamma^nN=\Gamma^n
$$
implies $\Gamma^n$ is $p$-divisible, so the fact that $N$ is a central extension
$$
1\rightarrow \Gamma^n\rightarrow N\rightarrow N/\Gamma^n\rightarrow 1
$$
concludes the proof.
\hfill $\blacksquare$

\begin{Lem} \label{ExtIsZero}
Suppose $N$ is a nilpotent group and $p$ is a prime.
The following conditions are equivalent:
\begin{itemize}
\item[a.] $\Ext(\ZZ/p^\infty,N)=0$,
\item[b.] $\Ext(\ZZ/p^\infty,N)$ is $p$-divisible,
\item[c.] $N$ is $p$-divisible.
\end{itemize}
\end{Lem}
\dokaz
(a)$\implies$(b) is obvious.
For (c)$\implies$(a) let $N$ be $p$-divisible. Then so is its abelianization and
Proposition 3 of \cite{CD2} implies $E=0$. 
\par
(b)$\implies$(c)
If $N$ is not $p$-divisible neither is its abelianization by  \ref{AbDivisibleImpliesGIsDivisible}.
Therefore, by Proposition 3 of \cite{CD2} 
${\mathrm{Ext}}(\mathbb{Z}_{p^{\infty}},{\mathrm{Ab}}N)$ is not $p$-divisible.
Then the six-term exact sequence of $Hom$ and $Ext$ [\cite{BK}, p.170] implies that
$\Ext(\ZZ/p^\infty,N)$ is not $p$-divisible.
\hfill $\blacksquare$

\begin{Lem} \label{HomIsZero}
Suppose $G$ is a nilpotent group and $p$ is a prime.
The following conditions are equivalent:
\begin{itemize}
\item[a.] $\Hom(\ZZ/p^\infty,G)=0$,
\item[b.] $\Hom(\ZZ/p^\infty,G)$ is $p$-divisible,
\item[c.] $\Hom(\ZZ/p^\infty,G)\otimes \ZZ/p^\infty=0$
\item[d.] $G$ does not contain $\ZZ/p^\infty$.
\end{itemize}
\end{Lem}
\dokaz 
Note that albeit Bousfield and Kan \cite{BK} defined $Hom$ as a space, they showed that it is
also the set of the respective homomorphisms.

(a)$\implies$(b) and (b)$\implies$(c) are obvious.
\par (c)$\implies$(b). Notice the $p$-torsion of $\Hom(\ZZ/p^\infty,G)$ is trivial.
Indeed, if $i:\ZZ/p^\infty\to G$ and $i^p=1$, then for any $a\in \ZZ/p^\infty$
we find $b\in \ZZ/p^\infty$ satisfying $b^p=a$.
Now, $i(a)=i(b^p)=i^p(b)=1$. If an Abelian group $A$ has no $p$-torsion
and $A\otimes \ZZ/p^\infty=0$, then $A$ is $p$-divisible.
\par 
(b)$\implies$(d) Suppose $i:\ZZ/p^\infty\to G$ is a monomorphism.
Given $a\in \ZZ/p^\infty$ find $k\ge 1$ such that $a^{p^k}=1$ and choose
$\phi:\ZZ/p^\infty\to G$ so that $i=\phi^{p^k}$.
Now $i(a)=\phi^{p^k}(a)=(\phi(a))^{p^k}=\phi(a^{p^k})=\phi(1)=1$, a contradiction.
\par (d)$\implies$(a). Given a non-trivial $i:\ZZ/p^\infty\to G$ its image
is a direct sum of copies of $\ZZ/p^\infty$, a contradiction.
\hfill $\blacksquare$

\begin{Lem} \label{PiOneIsPLocal}
Suppose $L$ is a nilpotent CW complex and $p$ is a prime.
If $H_1(L;\ZZ/p^\infty)=H_2(L;\ZZ/p^\infty)=0$,
then $\pi_1(L)$ is $\bar p$-local.
\end{Lem}
\dokaz In view of $H_2(L;\ZZ/p^\infty)=0$, $H_1(L)$ has trivial $p$-torsion and $H_1(L;\ZZ/p^\infty)=0$
implies $H_1(L)$ is $p$-divisible. So is $\pi_1(L)$ (see \ref{AbDivisibleImpliesGIsDivisible}).
Consider the $p$-completion $\hat L_p$ of $L$. As $\pi_1(\hat L_p)=\Ext(\ZZ/p^\infty,\pi_1(L))=0$
and $H_2(\hat L_p;\ZZ/p^\infty)=H_2(L;\ZZ/p^\infty)=0$ one gets
$\pi_2(\hat L_p)\otimes \ZZ/p^\infty=0$ by the Hurewicz Theorem.
The exact sequence $$0\to \Ext(\ZZ/p^\infty,\pi_2(L))\to \pi_2(\hat L_p)\to \Hom(\ZZ/p^\infty,\pi_1(L))\to 0$$
implies $\Hom(\ZZ/p^\infty,\pi_1(L))\otimes \ZZ/p^\infty=0$.
By \ref{HomIsZero} $\pi_1(L)$ is $\bar p$-local.
\hfill $\blacksquare$


\begin{thebibliography}{99}

\bibitem{Bousfield}
A.K.Bousfield, {\em Localization and periodicity in unstable
homotopy theory}, J.~Amer.~Math.~Soc.~{\bf 7} (1994), no.~4,
831--873.

\bibitem{BK}
A.K.Bousfield and D.M.Kan, {\em Homotopy limits completions and localizations}, 
Springer Lecture Notes in Math., Vol.304 (2nd corrected printing),
Springer-Verlag, Berlin-Heidelberg-New York, 1987.

\bibitem{CDMV} M.Cencelj, J.Dydak, A.Mitra,  and A.Vavpeti\v c,
{\em Bockstein Theorem for nilpotent groups}, preprint

\bibitem{CD2}
M.Cencelj and A.N.Dranishnikov,
{\em Extension of maps to nilpotent spaces II},
Topology Appl. 124 (2002), no. 1, 77--83.

 \bibitem{CD3}
M.Cencelj and A.N.Dranishnikov,
{\em Extension of maps to nilpotent spaces III},
Topology Appl. 153 (2005) 208--212.

 \bibitem{D-T}
 A.Dold and R.Thom,
 {\em Quasifaserungen und Unendliche Symmetrische Produkte},
 Annals of Math. 67 (1958),239--281.

\bibitem{DR0}
A.N.Dranishnikov, {\em Extension of mappings into CW complexes}, 
Mat. USSR Sbornik, 74 (1993), 47--56.


\bibitem{DR4}
A.N.Dranishnikov, {\em Cohomological Dimension Theory of Compact
Metric Spaces}, Topology Atlas (1999).

\bibitem{DrDy2}
 A. N. Dranishnikov and J. Dydak, {\em Extension dimension and extension types}, Proc.
Steklov Math. Inst. 212 (1996), 55--88.

\bibitem{DrDy1}
A.Dranishnikov and J.Dydak, {\em Extension theory of separable metrizable spaces
 with applications to dimension theory},  Transactions of the American Math.Soc. 353 (2000),
133--156.

 \bibitem{Dy1}
 J.Dydak,
{\em Cohomological dimension and metrizable spaces},
 Transactions of the Amer.Math.Soc. 337 (1993),219--234.

 \bibitem{Dy2}
 J.Dydak,
{\em Cohomological dimension and metrizable spaces},
 Transactions of the Amer.Math.Soc. 348 (1996), 1647--1661.
 
\bibitem{GoerssJardine}
P.G.Goerss and J.F.Jardine, {\em Simplicial Homotopy Theory}, Birkhauser,
Basel-Boston-Berlin, 1999.

\bibitem{Hatcher}
A.Hatcher, {\em Algebraic Topology}, Cambridge University Press,
Cambridge, 2002.


\bibitem{HMR}
P.~Hilton, G.~Mislin, J.~Roitberg, {\em Localization of Nilpotent
groups and spaces}, North-Holland Publishing Co.~,
Amsterdam-Oxford; American Elsevier Publishing Co.~, Inc.~, New
York, 1975.

\bibitem{Ku}
  W.I.Kuzminov, {\em Homological dimension theory},
  Russian Math. Surveys 23 (1968),  1--45.
  
\bibitem{Mac}
  J. MacCleary, {\em A User's Guide to Spectral Sequence},
  Cambridge studies in advanced mathematics, second edition, 2  

\bibitem{M-S}
S.Marde\v{s}i\'{c} and J.Segal, {\em Shape theory},
North-Holland Publ.Co. (1982).

\bibitem{Nei}
J.A.Neisendorfer, {\em Primary homotopy theory},
Memoirs.Amer.Math.Soc 232 (1980).

\bibitem{Robinson}
D.~J.~S.~Robinson, {\em A course in the theory of groups},
Springer-Verlag, New York, 1993.

\bibitem{Sullivan}
D.Sullivan, {\em Genetics of homotopy theory and the Adams
conjecture}, Ann. of Math. (2) {100} (1974), 1--79.

\bibitem{War}
R.B.Warfield, {\em Nilpotent groups}, 
Springer Lecture Notes in Math., Vol.513,
Springer-Verlag, Berlin-Heidelberg-New York, 1976.


\end{thebibliography}
\end{document}